\documentclass[11pt,reqno]{article}

\title {The speed of biased random walk on percolation clusters}

\author{Noam Berger\thanks{Research supported by Microsoft Research
 graduate fellowship.} \\ U.C. Berkeley \\ {\small noam@stat.berkeley.edu}
\and Nina Gantert\thanks{Research partially supported by the DFG under
 grant SPP 1033.} \\ University of Karlsruhe \\ {\small N.Gantert@math.uni-karlsruhe.de }
\and Yuval Peres\thanks{Research partially supported by NSF grant \#DMS-0104073
and by a Miller Professorship at UC Berkeley.} \\ U.C. Berkeley \\
 {\small peres@stat.berkeley.edu}
}

\usepackage{latexsym}
\usepackage{amsfonts}
\usepackage{amsmath}
\usepackage{amsthm}
\usepackage{epsfig}

\newtheorem{claim}{Claim}
\newtheorem{fact}[claim]{Fact}
\newtheorem{thm}{Theorem}
\newtheorem{defn}{Definition}
\newtheorem{lemma}{Lemma}
\newtheorem{cor}[thm]{Corollary}
\newtheorem{prop}[thm]{Proposition}

\newtheorem{remark}{Remark}

\newtheorem{conjecture}{Conjecture}

\newcommand{\hetzi}{\frac{1}{2}}

\newcommand{\Z}{{\mathbb Z}}

\newcommand{\olle}{{H\"{a}ggstr\"{o}m }}
\newcommand{\ve}{\text{ and }}

\renewcommand{\H}{{\cal H}}
\newcommand{\raz}{{\widetilde Z_n}}

\newcommand{\lattice}{\Z^2}
\newcommand{\sublattice}{\Z}
\newcommand{\PP}{{\mathbb P}}
\newcommand{\tildd}{{\widetilde D}}
\newcommand{\RT}{{\mbox{RT}}}
\newcommand{\OK}{{\mbox{OK}}}
\newcommand{\linezero}{{\cal L}}

\newcommand{\lP}{{\widehat{ P}_p}}

\begin{document}

\maketitle
\noindent

\begin{abstract}
We consider biased random walk on 
supercritical percolation clusters in $\lattice$. 
We show that the random walk is transient and that there are two speed 
regimes: If the bias is large enough, the random walk 
has speed zero, while if the
bias is small enough, the speed of the random walk is positive. 
\end{abstract}



\section{Introduction}
The following model is considered in the physics literature as a model for 
transport in an inhomogeneous medium. Let $p \in (p_c, 1)$,
where $p_c= \hetzi$ is the critical probability 
for bond percolation
on $\lattice$.
We perform i.i.d. bond percolation with parameter $p$ on $\lattice$. For
convenience, we always condition on the event that the origin belongs to the
infinite cluster. The 
corresponding measure on percolation configurations will be denoted 
by $P_p^*$. Let $\beta>1$. Consider the random walk starting at the
origin with transition probabilities defined as follows.
Let $Z_n=(X_n,Y_n)$ be the location at time $n$. 
Let $l_n$ be the
number of neighbors $Z_n$ has in the infinite cluster.
If $\raz=(X_n+1,Y_n)$ is one of these neighbors, then $Z_{n+1}=\raz$
with probability
\begin{equation*}
\frac{\beta}{\beta+l_n-1}
\end{equation*}
and $Z_{n+1}$ is any of the other neighbors with probability
\begin{equation*}
\frac{1}{\beta+l_n-1}\, .
\end{equation*}
If $\raz$ is not a neighbor of $Z_n$ (i.e. the edge $(Z_n,\raz)$ is closed)
then $Z_{n+1}$ is chosen among the neighbors of $Z_n$ with equal probabilities.
This is a random walk with bias to the right, where the strength of the bias is given by the parameter $\beta$. To our best knowledge, the first authors who 
considered this model are M. Barma and D. Dhar in \cite{dharbarma}.

Let $\omega$ be the percolation
configuration.
We write $P^\beta _\omega$ for the conditional law
of the random walk
given $\omega$, and $\PP^{\beta, *} $ for the joint distribution of 
$\left(\omega,(Z_n)_{n=1,2, \ldots}\right)$.
$\PP^{\beta,*} $ restricted to $(Z_n)_{n=1,2, \ldots}$ is the law of the walk
 averaged over the realizations of the 
percolation configuration. 

Our main result is the following theorem, which proves part of the
predictions of \cite{dharbarma}.
\begin{thm}\label{main}
For every $ p \in (p_c, 1)$,  
there exist $1<\beta_\ell\leq \beta_u<\infty$ such that
if $1<\beta<\beta_\ell$ then
\begin{equation*}
\lim_{n\to\infty}\frac{X_n}{n}>0\quad \PP^{\beta,*}\mbox{--a.s.}
\end{equation*}
and if $\beta>\beta_u$ then
\begin{equation*}
\lim_{n\to\infty}\frac{X_n}{n}=0\quad \PP^{\beta,*}\mbox{--a.s.}
\end{equation*}
\end{thm}

The following conjecture goes back to \cite{dharbarma}.
\begin{conjecture}\label{conj1}
The statements of Theorem \ref{main} hold with
$\beta_{\rm crit}: =\beta_\ell=\beta_u.$
\end{conjecture}

While there is a large physics literature on this model, 
as, for instance, \cite{dharbarma, dhar, stauffer}, there are 
few mathematical results. The biased random walk on the percolation
cluster is a random walk in a random environment on $\lattice$. There has been 
remarkable recent progress on laws of large numbers for 
random walk in dependent random environments, see 
\cite{cometszei, rassoul, shen}. However, in all of these papers, there are 
boundedness assumptions on the transition probabilities which are 
violated in our case.

In contrast to the biased case, simple random walks on percolation clusters were investigated in the
probability literature for some time. The first work on the subject
was done By Grimmett, Kesten and Zhang \cite{GKZ}, 
where they proved that
simple random walk on supercritical percolation clusters in $\Z^d$ is 
transient for $d\geq 3$. Other
papers include \cite{BPP}, \cite{EO}, \cite{hohe}, \cite{ABBP}, \cite{bemo}.

In order to prove 
that there is a 
positive speed regime, we first assume that $p$ is close enough to $1$ and 
show 
the following.
\begin{prop}\label{mainlargep}
For every $ p$ close enough to $1$,  
there exists $\beta_\ell > 1$ such that if $\beta<\beta_\ell$ then
\begin{equation}
\lim_{n\to\infty}\frac{X_n}{n}>0\quad \PP^{\beta,*}\mbox{--a.s.}
\end{equation}
\end{prop}
The paper is organized as follows. Sections \ref{sec:existence}, \ref{sec:traps}, \ref{observe}, \ref{sec:sqrt}, 
\ref{regenerations}, and \ref{proflrgp} are devoted to the proof of 
Proposition \ref{mainlargep}. 
Using 
renormalization arguments we show in Section \ref{renorm} that the 
statement of
Proposition \ref{mainlargep} holds for every $p> p_c$. In Section
\ref{speedzero}, we define $\beta_u$ and show that for $p > p_c$ and 
$\beta> \beta_u$, the speed is zero. In fact, our $\beta_u$ is the predicted 
value
of $\beta_{\rm crit}$ in Conjecture \ref{conj1}, 
see \cite{dharbarma}. The proofs in this section carry over to the 
multidimensional case, i.e. to biased random walks on supercritical percolation clusters in
$\Z^d$, $d\geq 2$.

While finishing this paper, we learned that A.\ S.\ Sznitman has 
independently obtained results similar to ours. In \cite{sznitman}, he 
investigates biased random walks on supercritical percolation clusters 
in $\Z^d$, $d\geq 2$, where the transition probabilities correspond to 
weights given by scalar products with a direction vector. He shows the
analogue of Theorem \ref{main} and obtains a CLT in the positive speed regime.
While both \cite{sznitman} and this paper use a regeneration structure to derive the main results, the techniques of the two papers are quite different.
Sznitman uses very precise information about the random walk and its 
analytical properties, while our approach uses more detailed information about the percolation cluster.

{\bf Acknowledgment} 
We thank Alan Hammond and Manjunath Krishnapur for useful
discussions. We also thank the referee for
helpful comments on the first version of the paper. 


\section{A positivity criterion for the speed}\label{sec:existence}
\begin{figure}
\center{\epsfig{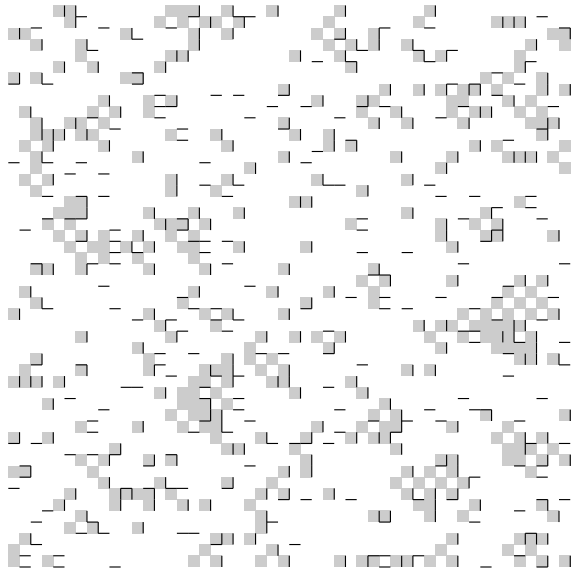}}
\caption{Traps in the percolation cluster, $p=0.9$. Vertices are represented by squares. Bad
vertices are shaded. Marked are the closed dual bonds}\label{fig:traps}
\end{figure}         

In order to simplify the arguments, we 
will without loss of generality condition, throughout the proof, on the 
event that 
the origin is in the infinite cluster on the left half--plane, i.e. on
the event that there is an infinite cluster on $\{(x,y): x\leq 0 \}$
and that the origin is in this infinite cluster. 
This event has positive probability: in fact, due to the results of 
\cite{BGN}, the probability that there is an 
infinite cluster on the left half--plane equals $1$ whenever 
$p >p_c$ (see also \cite{grimmett}). 
Denote the corresponding probability measure on percolation 
configurations by $\lP$, and the resulting joint law of 
$\left(\omega,(Z_n)_{n=1,2, \ldots}\right)$ by $\PP^\beta$.
\begin{remark}
Assume that the origin is not in the 
infinite cluster on 
the left half--plane. Take an arbitrary vertex $z$ which {\bf is} in 
the infinite 
cluster to its left. 
Then there is a finite open path $\Gamma$ connecting $z$ to the origin. 
If the statements in Theorem \ref{main} hold almost surely for 
the random walk starting from $z$, then they also hold almost surely
for the random walk starting from the origin, since
starting from $z$, we have a 
positive probability to go to the origin.

\end{remark}

We give a criterion which will later be used to show that the speed is 
strictly positive for $\beta$ small enough.
We will prove in Lemma \ref{transtoright} and Lemma \ref{rentranstoright} that
\begin{equation}\label{toright}
\lim_{n\to\infty}X_n=\infty \quad\PP^\beta\mbox{--a.s.}
\end{equation}
We call $n>0$ a {\bf fresh epoch} if $X_n > X_k$ for all $k< n$ and
we call $n$ a {\bf regeneration epoch } if, in addition, 
$X_k > X_{n}$ for all $k> n$. Let the regeneration epochs be $0=R_0 < R_1 < R_2 < \ldots$.
Exactly as in \cite{LPP}, one shows that there are, $\PP^\beta$--a.s., infinitely many regeneration epochs and that the time differences $(R_{i+1} - 
R_i)_{i=1, 2, 3, \ldots}$ and the increments between regeneration epochs
$(X_{R_{i+1}} - X_{R_i})_{i=1, 2, 3, \ldots}$ are i.i.d. sequences 
under $\PP^\beta$.
Standard arguments then imply that if
$E^\beta(R_{2} - R_1) < \infty$, then
\begin{equation}\label{speedf}
\lim_{n\to \infty} \frac{X_n}{n} = \frac{E^\beta(X_{R_{2}} - X_{R_1})}
{E^\beta(R_{2} - R_1)}>0 \quad \PP^\beta\mbox{--a.s.}
\end{equation}

\section{An exponential bound on the size of traps} \label{sec:traps}
\begin{figure}
\center{\epsfig{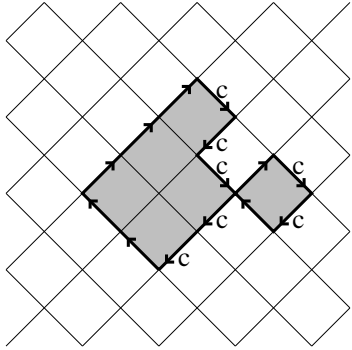}}
\caption{The unique contour around an even  trap in the even
  lattice.
Vertices of the even lattice are represented as squares. 
Marked are the dual bonds.}\label{contour}
\end{figure}         

We use the following decomposition of the percolation cluster into
good and bad points. The definition of a good point might
seem artificial at first sight, but the results of Sections
\ref{observe} and \ref{sec:sqrt} will clarify the choice of this definition.
\begin{defn}\label{def:goodpnt}
A point $z=(x,y)\in\lattice$ is {\bf good} if there exists an infinite path
$\{z_0=z,z_1=(x_1,y_1),z_2=(x_2,y_2),...\}$  
such that for $k=1,2,3, \ldots $,
\\(A) $|y_k-y_{k-1}|=1$ and $x_k-x_{k-1}=1$.
\\(B) The edges
$\{(x_{k-1},y_{k-1}),(x_{k},y_{k-1})\}$
and $\{(x_{k}, y_{k-1}) , (x_{k},y_{k})\}$ are open.
\end{defn}
Denote the infinite cluster by $I$ and the set of good vertices by $J$.
A vertex $z$ is {\bf bad} if $z \in I$ and $z$ is not good. Connected components of $I \setminus J$
will be called {\bf traps} (see Figure \ref{fig:traps}).
For a vertex $v$, let $C(v)$ be the trap containing $v$. 
($C(v)$ is empty if $v$ is a good point.) The length of
a trap $T$ is
\begin{equation*}
L(T)=\sup\{|x_1-x_2|:\exists y_1,y_2 \text{ such that }(x_1,y_1)\in T
\ve (x_2,y_2)\in T \}
\end{equation*}
and the width is
\begin{equation*}
W(T)=\sup\{|y_1-y_2|:\exists x_1,x_2 \text{ such that }(x_1,y_1)\in T
\ve (x_2,y_2)\in T \}
\end{equation*}
If $T$ is empty, then we take $L(T)=W(T)=0$. For convenience, we
will use the notation $L(v)$ for $L(C(v))$ and $W(v)$ for
$W(C(v))$.
\begin{lemma}\label{trapsize}
For every $p$ close enough to $1$,  there exists $\alpha=\alpha(p)<1$ such that
$\lP(L(0)\geq n)\leq \alpha^n$ and $\lP(W(0)\geq n)\leq \alpha^n$ for
every $n$. Further, $\lim_{p\to 1}\alpha(p) = 0$.
\end{lemma}
\begin{proof}
Call two vertices {\bf even--connected} if $\|u-w\|_1=2$. That is, $(x,y)$
and $(x',y')$ are even--connected if either $|x-x'|=|y-y'|=1$ or
$\left(|x-x'|,|y-y'|\right)=(0,2)$ or 
$\left(|x-x'|,|y-y'|\right)=(2,0)$.
We define the
{\bf even trap} \label{pageeventrap}
$C_e(v)$ of a point $v$ as the
even--connected component of bad points $v'$ with
\[\|v'\|_1\equiv \|v\|_1 \mod 2,\]
containing $v$. 
 In particular, all points $v'$ in $C(v)$ with
\[\|v'\|_1\equiv \|v\|_1 \mod 2\]
are also in the even trap of
$v$ (but the even trap may contain additional points not in $C(v)$). The following is obvious.
\begin{fact}\label{eventrap}
Every vertex in $C(0)$ is either an element of $C_e(0)$ or a
neighbor of a vertex in $C_e(0)$. In
particular, $L(0)\leq L(C_e(0))+2$ and $W(0)\leq W(C_e(0))+2$.
\end{fact}
Thanks to Fact \ref{eventrap}, we only need to give exponential
bounds to $L(C_e(0))$ and to $W(C_e(0))$. Consider the following
percolation model on the even lattice (i.e. the lattice whose vertices
are $\{v\in \lattice :\,  \|v\|_1 \mbox{ is even}\}$ and which has an
undirected edge between every $(x,y)$ and $(x',y')$ such that $|x-x'|=|y-y'|=1$):
The bond between $(x,y)$ and $(x+1,y\pm 1)$ is open if and
only if in the original model the edges $\{(x,y),(x+1,y)\}$ and
$\{(x+1,y),(x+1,y\pm 1)\}$ are open. 
This is a model of dependent oriented percolation, and we denote
the corresponding probability measure by $P_{p,\rm oriented}$.

Let $p'$ be close to $1$. By the results in \cite{lss}, there exists $p<1$ such that
$P_{p,\rm oriented}$ dominates i.i.d. bond percolation 
with parameter $p'$ on the even lattice. Consider $C_e(0)$ in the even lattice.
Let the outer boundary of a 
set of vertices be the set of all edges which have one vertex in the 
set and one in the complement. The outer boundary can be identified with 
a contour in the dual
lattice (see Figure \ref{contour}). Hence, the number of 
outer boundaries of size $n$ is bounded by $\exp(O(n))$ (each contour, which 
needs not to be simply connected, can be identified with a random walk path).
By an argument similar to that of
\cite{duret} p. 1026, at least half of the edges in the 
outer boundary of $C_e(0)$
are closed (in Figure \ref{contour}, these are the boundary edges marked with
a ``C'').
Therefore, if $p'$ is close enough to $1$, 
$L(C_e(0))$ and $W(C_e(0))$ have the desired exponential tail
with respect to i.i.d. bond percolation
with parameter $p'$ on the even lattice, hence 
also with respect to $P_{p,\rm oriented}$.
\end{proof}

\section{Bound for back--stepping from a good vertex}\label{observe}
The following simple observation is essential to the proof. 
Let ${\cal H}(n)$ be the $\sigma$--field generated by the history
of the random walk until time $n$, i.\ e.\ 
${\cal H}(n) = \sigma(\{Z_0 = 0, Z_1, Z_2, \ldots ,Z_n\})$.  
Let $P^\beta_{\omega, {\cal H}(n)}$ be the conditional 
distribution of
$P^\beta_{\omega}$, given ${\cal H}(n)$, and $\PP^\beta_{{\cal H}(n)}$ be 
the conditional 
distribution of
$\PP^\beta$, given ${\cal H}(n)$.
Define $\tau_n(X)=\min\{i>n:X_n=X\}$. In order not to overload the
notations, in many places throughout the section we chose to omit the
integer brackets, e.g. we write $\ell/3$ instead of 
$\left[\ell/3\right]$.
\begin{lemma}\label{preconduc}
There exists $D'=D'(\beta)$ such that for
every $\ell=1,2,3, \ldots$ and for every configuration $\omega$ such that $z=(x,y)$ is
a good point, 
\begin{equation*}
P^\beta_{\omega, {\cal H}(n)}(\tau_n(X_n-\ell)\leq\tau_n(X_n+\ell/3)|Z_n=z )<
D'\beta^{-\ell/3}.
\end{equation*}
\end{lemma}
\begin{proof}
The transition probabilities can be described with the following 
electrical network: Give a weight to each open edge  $e$:
if $e=\{(x,y),(x+1,y)\}$ then $e$ has weight
$w(e)= \beta^{x+1}$, and if $e=\{(x,y),(x,y\pm 1)\}$ then $e$ has weight
$w(e)= \beta^x$. If $e$ is closed, then its weight is $0$. The random walk 
$(Z_n)$ has transition probabilities proportional to the weights of the 
edges from a vertex. For background on the description of reversible Markov chains as electrical networks, we refer to \cite{doylesnell} and to 
\cite{lyons}.

The following fact is well known, but for the convenience of the
reader we will recall its proof.
\begin{fact}\label{fact:conduc}
Let $G$ be a finite electrical network, and let $A$ and $B$ be disjoint
sets of vertices in $G$. Let $z$ be a vertex in $G$, and let
$\tau(z\to A)$ (resp. $\tau(z\to B)$) be the hitting time of $A$
(resp. $B$) for a walk starting at $z$. Let $C_{z,A}$
(resp. $C_{z,B}$) be the effective conductance between $z$ and $A$
(resp.  $B$). Then,
\begin{equation}\label{hitbnd}
P\left(\tau(z\to B)<\tau(z\to A)\right)\leq\frac{C_{z,B}}{C_{z,A}}.
\end{equation}
\end{fact}
\begin{proof}
Let $\pi(z)$ be the sum of the weights of all edges
containing $z$. Let $u_j$ be the location of the walker at time $j$.
Let $k_i$ be the $i$--th time the walk returns to $z$
(i.e. $k_0=0$, and $k_{i+1}=\tau_{k_i}(z)$). We call the interval
$[k_{i-i},k_i-1]$ the $i$--th excursion.
For a set $D\subseteq G$, let $V(i,D)$ be the event that the walker
visits $D$ during the $i$--th excursion.
Then, for every $i$,
\begin{equation}\label{exvis}
P(V(i,D))=\frac{C_{z,D}}{\pi(z)}.
\end{equation}
(see e.g. equation (2.4) of \cite{lyons}).  By (\ref{exvis}), for every $i$,
\begin{equation*}
P\left(V(i,B)|V(i,A\cup B)\right) = \frac{C_{z,B}}{C_{z,A\cup B}}.
\end{equation*}
In particular, decomposing the sequence of excursions according to 
the first visit to $A\cup B$ and using the fact that the excursions are i.i.d.,
we get
\begin{equation*}
P\left(\tau(z\to B)<\tau(z\to A)_{\text{ }}\right)
\leq
\frac{C_{z,B}}{C_{z,A\cup B}}
\leq\frac{C_{z,B}}{C_{z,A}}.
\end{equation*}
\end{proof}

Consider the box $B=[x-\ell,x+\ell/3]\times
[y-\beta^{2\ell/3},y+\beta^{2\ell/3}]$. In view of Fact \ref{fact:conduc}, we
need to estimate the effective conductances between $z$ and the face
$B^+=\{x+\ell/3\}\times[y-\beta^{2\ell/3},y+\beta^{2\ell/3}]$ 
and between $z$ and the rest of the boundary of the rectangle.
\begin{enumerate}
\item\label{bplus}
$C_{z,B^+}$ is bounded from below by the conductance of the good path
from $z$ to $B^+$, which is at least $D_1\beta^x$ for some $D_1=D_1(\beta)$.
\item\label{bminus}
Consider $B^-=\{x-\ell\}
\times[y-\beta^{2\ell/3},y+\beta^{2\ell/3}]$. The conductance $C_{z,B^-}$
is bounded from above by the sum of the weights $w(u,u+(1,0))$
for $u\in B^-$. But for every such $u$, $w(u,u+(1,0))\leq \beta^{x-\ell+1}$
(with inequality because the weight is zero if the edge is closed),
and there are $2\beta^{2\ell/3}$ such edges. Therefore
$C_{z,B^-}\leq D_2\beta^x\cdot \beta^{-\ell/3}$ for some $D_2=D_2(\beta)$.
\item\label{bstar}
Consider $B^*_1=[x-\ell,x+\ell/3]\times\{y-\beta^{2\ell/3}\}$ and
$B^*_2=[x-\ell,x+\ell/3]\times\{y+\beta^{2\ell/3}\}$. By Nash--Williams'
inequality (equation (2.15) on page 38 of \cite{lyons}),  
\begin{equation*}
C_{z,B^*_j}\leq\beta^{-2\ell/3}\sum_{i=x-\ell}^{x-1+\ell/3}\beta^{i+1}\leq
D_3\beta^x\cdot\beta^{-\ell/3}
\end{equation*}
for some $D_3=D_3(\beta)$.
\end{enumerate}
From \ref{bplus}., \ref{bminus}. and \ref{bstar}. we see, using 
(\ref{hitbnd}), that the
probability to exit $B$ not through $B^+$ is at most $O(\beta^{-\ell/3})$.
\end{proof}

The following lemma gives a bound for the probability of back--stepping from 
a good point at a fresh epoch.
Recall that $n>0$ is a {\bf fresh epoch} if $X_n > X_k$ for all $k< n$.
\begin{lemma}\label{conduc}
Assume that $p$ is close enough to $1$.  Let $G(z)$ be the event
that $z$ is a good point and let $F(n)$ be the event that $n$ is a
fresh epoch. 
Then there exists $K=K(\beta,p)$ such that for every $\ell=1,2,\ldots$,
\begin{equation*}
\PP^\beta_{\H(n)}(\hbox {\rm there is an } m\geq n 
\hbox{ \rm such that } X_m\leq
x-\ell)\left|Z_n=z, F(n), G(z)\right.)\leq  
K\beta^{-\sqrt{\ell}/K},\, \, \PP^\beta\mbox{--a.s.}
\end{equation*}
\end{lemma}
To prove Lemma \ref{conduc}, we will use the following lemma:
\begin{lemma}\label{tromconduc}
In the notations of Lemma \ref{preconduc}, let $\tau'_n(X)$ be the
first fresh epoch, later than $n$, such that the random walk hits a
good point whose
first coordinate is larger or equal to $X$. Then, there exists a
constant $D=D(\beta,p)$ such that for every $\ell=1,2,\ldots$,
\begin{equation*}
\PP^\beta_{\H(n)}\left(\tau_n(X_n-\ell)<\tau'_n(X_n+\ell/6)\left|Z_n=z, G(z),
    \max_{0\leq i\leq n}X_i <X_n+ \sqrt{\ell}
\right.\right)
\leq D\beta^{-\sqrt{\ell}/D},\,\, \PP^\beta\mbox{--a.s.}
\end{equation*}
In particular, 
\begin{equation}\label{eq:tromconduc}
\PP^\beta_{\H(n)}(\tau_n(X_n-\ell)<\tau'_n(X_n+\ell/6)\left|Z_n=z, F(n), 
G(z)\right.)
\leq D\beta^{-\sqrt{\ell}/D},\quad \PP^\beta\mbox{--a.s.}
\end{equation}
\end{lemma}
\begin{proof}
For $i= 1,\ldots,[\sqrt{\ell}/6]$, let
$t_i=\tau_n(X_n+i\sqrt{\ell})$. For convenience, if
$t_i=\infty$ then we say that $Z_{t_i}=\infty$ and $Z_{t_i}$ is 
not a good point.
We define the {\bf right hand trap} (resp. {\bf right hand even trap}) of 
a bad point $z=(x,y)$ to be the connected
component (resp. even connected component) of bad points $z'=(x',y')$
such that $x'\geq x$, containing $z$. The right hand even trap 
of a point $z$ will be
denoted by $\RT(z)$.
Let $L(\RT(z))$ be the length of the (right hand even) trap
$\RT(z)$. If $z$ is a good point then we say that $L(\RT(z))=0$.
\begin{claim}\label{lengthrt}
For $z=(x,y)$, let $\omega_l(z)$ be the configuration of all 
edges to the left of the line $L_x= \{(x, \widetilde y)| \widetilde y 
\in \sublattice\}$, and
let $\omega_r(z)$ be the configuration of all edges to the right of $L_x$,
including the vertical edges on the line $L_x$.
For $\alpha$ as in Lemma \ref{trapsize}, and $k=1, 2, \ldots$,
\begin{equation}\label{trs}
\lP\left(L(\RT(z))\geq k\left|\mbox{ }\omega_l(z)\right.\right)
\leq \alpha^k, \quad \lP\mbox{--a.s.}
\end{equation}
In particular,
\begin{equation}\label{trs2}
\lP\left(G(z)\left|\mbox{ }\omega_l(z)\right.\right)\geq 1-\alpha,\quad 
\lP\mbox{--a.s.}
\end{equation}
\end{claim}
\begin{proof}
Since we condition on the origin being in the infinite cluster on the 
left half--plane, the event $\{L(\RT(z))\geq k\}$ is independent of 
${\omega_l(z)}$ and the claim follows from the proof of Lemma \ref{trapsize}.
\end{proof}
We want to estimate the probability of the following event: There exists
some $1\leq i\leq[\sqrt{\ell}/6]$ such that $t_i<\tau_n(X_n-\ell)$ and the
point $Z_{t_i}$ is good.
By (\ref{trs}), for every $i$, conditioned on $t_i<\infty$, 
\begin{equation}\label{strww}
\PP^\beta_{\H(t_i)}\left(L(\RT(Z_{t_i}))\geq \hetzi\sqrt{\ell}\right)\leq\alpha
^{\hetzi\sqrt{\ell}},\quad \PP^\beta\mbox{--a.s.}
\end{equation}
Using (\ref{trs2}), 
again conditioned on $t_i<\infty$, yields
\begin{equation}\label{gpww}
\PP^\beta_{\H(t_i)}\left(G(Z_{t_i})\left|L(\RT(Z_j)) < {\hetzi\sqrt{\ell}}
\mbox{ for all } 1\leq j<i
\right.\right)\geq 1-\alpha, \quad \PP^\beta\mbox{--a.s.}
\end{equation}
since we condition on an event which is measurable with respect to $\omega_l$.
The lemma now follows from Lemma \ref{preconduc}, 
(\ref{strww}) and (\ref{gpww}).
\end{proof}

\begin{proof}[Proof of Lemma \ref{conduc}]
Lemma \ref{conduc} now follows 
from (\ref{eq:tromconduc}) in Lemma \ref{tromconduc} by iterating.
\end{proof}

\section{An a priori bound}\label{sec:sqrt}
In this section we show an a priori bound for the distance the random walk
goes to the right.
\begin{lemma}\label{sqrt}
If $p$ is close enough to $1$, then for $\beta >1$  close enough to
$1$, there exists a constant $C$ such that for every $n$ large enough,
\begin{equation*}
\PP^\beta(X_n \leq Cn^\frac{1}{10})\leq n^{-2}.
\end{equation*}
\end{lemma}
In order to prove Lemma \ref{sqrt} we will give an estimate on the number 
of distinct sites visited by the random walk.

\begin{defn}
For a trap $T$, the {\bf size} of $T$ is $S(T)=L(T)+W(T)$.
\end{defn}

\begin{claim}\label{mix}
Let $T$ be a trap of size at most $s$, and let $z=(x,y)\in T$. Let 
\begin{eqnarray*}
\phi(s)&=& \beta^s \left(s^2 + \frac{2}{\beta -1}\right)\cdot
(3+\beta)\\ &\leq& C(\beta)\beta^{2s}.
\end{eqnarray*}
Then, for every $m$, and for every configuration $\omega$ with $z$ and $T$ 
as above,
\begin{equation*}
\PP^\beta_\omega\left( \#\{i:Z_i=z\}\geq m\right) \leq \left(1-\phi(s)^{-1}\right)^m.
\end{equation*}
In particular, if $z$ is a good point, then
\begin{equation*}
\PP^\beta_\omega\left( \#\{i:Z_i=z\}\geq m\right) \leq \left(1-\phi(1)^{-1}\right)^m.
\end{equation*}
\end{claim}

\begin{proof}
Recall the description of the transition probabilities with an electrical
network. By equation (2.3) of \cite{lyons}, starting at $z$, the probability of ever hitting $z$ again is
\begin{equation}\label{lyons}
1-\frac{C_{z,\infty}}{\pi(z)}
\end{equation}
where $\pi(z)$ is the sum of the weights of all edges containing $z$. Clearly,
\begin{equation}\label{piz}
\pi(z)\leq\beta^x(3+\beta).
\end{equation}
We need to bound $C_{z,\infty}$ from below. 
In order to do that we will
bound the resistance $R_{z,\infty}=1/C_{z,\infty}$ from above. For a good point
$z=(x,y)$, the resistance $R_{z,\infty}$ is bounded from above
by the resistance of the good path which is
\begin{equation}\label{rz0}
\frac{2 \beta^{-x}}{\beta -1}\ .
\end{equation}
If $z$ is in a trap $T$  of size at most $s$, let $z_0$ be a good
point on the boundary of $T$. Then,
\begin{equation}\label{sumr}
R_{z,\infty}\leq R_{z,z_0}+R_{z_0,\infty}.
\end{equation}
Let $q=(z,\ldots,z_0)$ be a path in $T$ from $z$ to $z_0$. Then, the
resistance of $q$ is bounded by the the product of the length of $q$ and the 
maximal
resistance of all bonds in $q$. Since $q$ is in $T$, its length is
bounded by $s^2$ and the maximal resistance of all bonds in $q$ is bounded by 
the maximal resistance of all bonds in $T$ which is at most $\beta^{s-x}$. 
Therefore,
\begin{equation*}
R_{z,z_0}\leq s^2\beta^{s-x}.
\end{equation*}
Further, $x_0\geq x-s$. Hence, using (\ref{rz0}) and (\ref{sumr}),
\begin{equation}\label{rz}
\frac{1}{C_{z,\infty}}=R_{z,\infty}\leq \beta^{-x}\cdot\left(s^2\beta^s + \frac{2\beta^s}{\beta - 1}\right)
\end{equation}
The claim now follows from (\ref{lyons}), (\ref{piz}) and (\ref{rz}).
\end{proof}

\begin{proof}[Proof of Lemma \ref{sqrt}]
Let $p$ be close enough to $1$ so that $P_p(S(0)\geq n)\leq \alpha^n$ for all $n$
large enough, with some $\alpha <1$.
Let $u$ be large enough so that
\begin{equation}\label{boundonalpha}
u\log \alpha<-4,
\end{equation}
and $\beta >1$ close enough to $1$ so that
\begin{equation}\label{boundonu}
u<\frac{1}{200\log\beta}\, .
\end{equation}
By (\ref{boundonu}), for every large enough $n$ and for
every $s\leq u\log n$,
\begin{equation}\label{choiceu}
\phi(s)\leq n^{1/10}.
\end{equation}
By the choice of $u$ the probability that there
exists a trap or an even trap of size bigger than $u\log n$ somewhere in the
square $[-n,n]\times[-n,n]$ is smaller than $\frac{1}{2}n^{-2}$.
We now condition on the event $A_1$ that there are no such traps.
Since
at times up to $n$ the random walk cannot leave the cube
$[-n,n]^2$, at any time before $n$ we are either
at a good vertex or in a trap of size at most $u\log n$.
\begin{claim}\label{numpoints}
Conditioned on $A_1$, with probability larger than $1-\exp(-\hetzi n^{1/5})$, 
the
random walk visits at least $n^{7/10}$ points up to time $n$.
\end{claim}
\begin{proof}
By (\ref{choiceu}) and by Claim \ref{mix}, for every
$z\in[-n,n]\times[-n,n]$, the probability that $z$ is visited more
than $n^{3/10}$ times is bounded by
\begin{equation*}
\left(1-n^{-1/10}\right)^{n^{3/10}}
\leq\exp\left(-n^{2/10}\right).
\end{equation*}
Therefore, the probability that any point in $[-n,n]\times[-n,n]$ is visited 
more than
$n^{3/10}$ times is bounded by
\begin{equation}\label{ubound}
4n^2\exp\left(-n^{2/10}\right)
\leq\exp\left(-\frac{1}{2} n^{1/5}\right)
\end{equation}
for $n$ large enough. But if no point is visited more than $n^{3/10}$
times, then at least $n^{7/10}$ points are visited.
\end{proof}
Let $B$ be the event that the
random walk visits at least $n^{7/10}$ points up to time $n$. 
\begin{claim}\label{diffx}
Conditioned on $B$, with probability at least $1-\exp(-\frac{1}{4}n^{1/5})$,
\begin{equation*}
\max\limits_{1\leq i \leq n} X_i- \min\limits_{1\leq i \leq n}X_i 
\geq n^{1/10}.
\end{equation*}
\end{claim}
\begin{proof}
Recall the Varopoulos--Carne bound for the $n$--step transition probabilities of a reversible Markov chain with reversible measure $\pi$ (see \cite{carne}):
\begin{equation}\label{vcarn}
P^n(a,b)\le 2 \sqrt{\pi(b)/\pi(a)} \exp\left(\frac{-d(a,b)^2}{2n}\right)
\end{equation}
where $d(a,b)$ is the (graph) -- distance between $a$ and $b$. 
Using (\ref{vcarn}) and a union bound, for every $1\leq i<j\leq n$, and all 
$\omega$
\begin{equation*}
P^\beta_\omega\left(X_i=X_j \ve |Y_i-Y_j|\geq n^{6/10}\right)\leq Cn^4\exp
\left(-\hetzi n^{1/5}\right)
\end{equation*}
where $C=C(\beta)$ is a constant. 
Taking the union over all possible pairs $i,j$,
\begin{eqnarray*}
P^\beta_\omega\left(\exists {i,j} \mbox{ such that } X_i=X_j \ve |Y_i-Y_j|\geq
  n^{6/10}\right)&\leq& Cn^6\exp\left(-\hetzi n^{1/5}\right)\\
&\leq& \exp\left(-\frac{1}{4} n^{1/5}\right)
\end{eqnarray*}
for $n$ large enough.
However, if $\max\limits_{1\leq i \leq n} X_i- \min\limits_{1\leq i \leq n}X_i
\leq n^{1/10}$ and at least $n^{7/10}$ points are visited, then there have 
to be $i$ and $j$ such that $X_i=X_j \ve |Y_i-Y_j|\geq n^{6/10}$.
\end{proof}
\begin{claim}\label{minx}
With probability at least $1-\exp(-n^{1/30})$, for every $1\leq i<j \leq
n$
\begin{equation*}
X_j-X_i \geq -n^{1/20}\, .
\end{equation*}
\end{claim}
\begin{proof}
For $z=(x,y)$ and $z'=(x',y')$, 
\begin{equation*}
\frac{\pi(z')}{\pi(z)}\leq C\beta^{x'-x},
\end{equation*}
where $C=C(\beta)$ is a constant. Fix $i<j$ and $z$ and $z'$ in $[-n,n]\times[-n,n]$ such that $x-x'>n^{1/20}.$
Then, again using (\ref{vcarn}),
\begin{equation*}
P^\beta_\omega(Z_i=z \ve Z_j=z') \leq 2\sqrt{\frac{\pi(z')}{\pi(z)}}\leq C\beta^{-n^{1/20}}
\end{equation*}
Summing over all of the possible values of $i$, $j$, $z$ and $z'$, we
get
\begin{eqnarray*}
P^\beta_\omega\left(\exists {i<j}\mbox{ such that }X_j-X_i<-n^{1/20}\right) &\leq&
C n^6\beta^{-n^{1/20}}\\ &\leq& \exp\left(-n^{1/30}\right)
\end{eqnarray*}
for $n$ large enough.
\end{proof}
Hence, with $\PP^\beta$--probability at least
$1-n^{-2}$, by Claim \ref{minx}
\[
\min\limits_{1\leq i \leq n}X_i
\geq-n^{1/20}.
\]
and, by Claim \ref{diffx},
\[
\max\limits_{1\leq i \leq n}X_i-\min\limits_{1\leq i \leq n}X_i \geq n^{1/10},
\]
hence
\[
\max\limits_{1\leq i \leq n}X_i\geq
n^{1/10}-n^{1/20},
\]
but, again due to Claim \ref{minx}, 
\[
X_n-\max\limits_{1\leq i \leq n}X_i \geq
-n^{1/20}.
\] 
Hence, with $\PP^\beta$--probability at least
$1-n^{-2}$, for $n$ large enough,
\begin{equation*}
X_n\geq n^{1/10}-2n^{1/20}\geq \hetzi n^{1/10}.
\end{equation*}
\end{proof}

\begin{lemma}\label{transtoright}
Let $p$ be close enough to $1$, and $\beta>1$. Then
\begin{equation}\label{gotorightas}
\lim_{n\to\infty}X_n=\infty \quad \PP^\beta\mbox{--a.s.}
\end{equation}
\end{lemma}

\begin{proof}
We prove the lemma by iterating Lemma \ref{tromconduc}.
Let $N>1$ be an arbitrary positive integer.
Let $T$ be the even trap containing the origin.
Let $\ell_0=2L(T)^2+N$, and let $\ell_{i+1} = 13\ell_i/12$ for every
$i=0,1,\ldots$. Let $\tau_g$ be the first time in which the walker is
in a good point. We define inductively the following times:
$t_0=\tau_g$,
$t_{i+1}=\tau'_{t_i}\left(X_{t_i}+ \ell_i/6\right)$.
Let
\begin{equation*}
A_i = 
\left\{\tau'_{t_i}\left(X_{t_i}+\ell_i/6\right)<
\tau_{t_i}(X_{t_i}-\ell_i)\right\}
\end{equation*}
Then by Lemma \ref{tromconduc}, for every $i$,
\begin{equation*}
\PP^\beta(A_i) \geq 1 - D\beta^{-\sqrt{\ell_i}/D}.
\end{equation*}
(The first formula in Lemma \ref{tromconduc} is needed since $t_0$ is not necessarily a fresh epoch). 
Therefore,
\begin{equation}\label{tlutbn}
\PP^\beta\left(\bigcap_{i=1}^\infty A_i\right) \geq 1 - 2C\beta^{-\sqrt{N}/D}
\end{equation}
for some $C=C(\beta)$. 
Note that $X_{t_i} - \ell_i \geq X_{t_{i-1}}- \frac{11}{12}\ell_{i-1}$. 
Hence, if $A_i$ occurs for every $i$, then $t_i<\infty$ for every
$i$, and 
\[
X_{s}> X_{t_i}-\ell_i\geq X_{t_0}-\ell_0+ \frac{1}{12}
\sum\limits_{j=1}^{i-1}\ell_j  
\geq X_{t_0} + C\ell_0 \left(\frac{13}{12}\right)^i
\]
for every $s>t_i$. In particular, if $A_i$ occurs for every $i$ then
(\ref{gotorightas}) holds. By (\ref{tlutbn}), the event in 
(\ref{gotorightas}) occurs with 
probability at least $1 -
2C\beta^{-\sqrt{N}/D}$ for every $N$. Therefore, 
it occurs a.s.
\end{proof}

\section{The environment after a regeneration}\label{regenerations}

\begin{figure}
\center{\epsfig{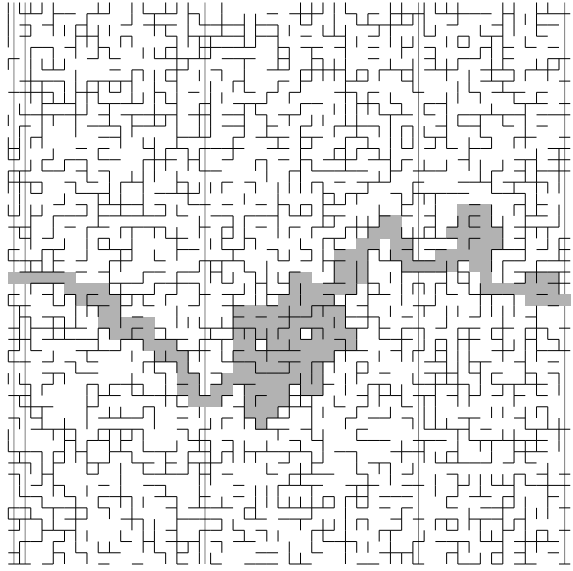}}
\caption{A sample path,
$p=0.55$ and $\beta=2.5$, gray lines are at regenerations.}
\label{fig:reg1}
\end{figure}

\begin{figure}
\center{\epsfig{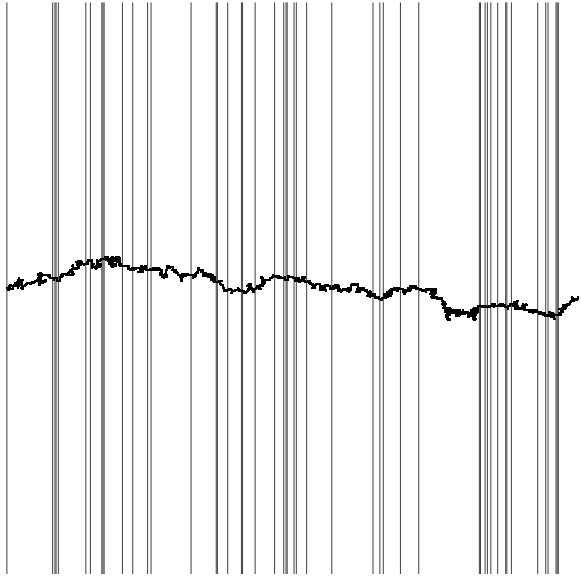}}
\caption{A sample path,
$p=0.65$ and $\beta=2.5$, gray lines are at regenerations.}
\label{fig:reg2}
\end{figure}

As the reader recalls from Section \ref{sec:existence}, we say that
$n>0$ is a {\bf fresh epoch} if $X_n > X_k$ for all $k< n$ and we say
that a fresh epoch $n$ is a {\bf regeneration epoch} or 
{\bf regeneration} if, in addition, 
$X_k > X_{n}$ for all $k> n$ (see Figures \ref{fig:reg1} and \ref{fig:reg2}).
In this section we consider the distribution of the percolation cluster to the right of $Z_n =z$, given that $n$ is a regeneration.

For every $z=(x,y)\in\lattice$, let $S_z$ be the environment to the
right of $z$, i.e. for every $\widetilde x>0$ and 
$\widetilde y\in\sublattice$, the edge
$S_z(\{(\widetilde x,\widetilde y),(\widetilde x,\widetilde y \pm 1)\})$ is open if and only if 
$\{(x+\widetilde x,y+\widetilde y),(x+\widetilde x,y+ \widetilde y \pm 1)\}$ 
is open and 
$S_z(\{(\widetilde x,\widetilde y),(\widetilde x+1,\widetilde y)\})$ is open 
if and only if 
$\{(x+\widetilde x,y+\widetilde y),(x+\widetilde x+1,y+\widetilde y)\}$ is
open.
For every time $n$ let $(F_n)$ be the future of the walk after time $n$,
i.e. $F_n(k)=Z_{n+k}-Z_{n},k=0,1,2, \ldots $. Let $\mu= \mu^\beta$ be the distribution of
$\{(F_0),S_0\}$ under $\PP^\beta$, conditioned on the event $\{X_i\geq
1\,  \forall i\geq 1 \}$. This is well defined because $\PP^\beta(X_i\geq
1\,  \forall i\geq 1) >0$.

\begin{lemma}\label{preregen}
Let $R_n$ be the $n$--th regeneration. Then, for all $n$, the law of
$\{(F_{R_n}),S_{Z_{R_n}}\}$ is $\mu$.
\end{lemma}
Lemma \ref{preregen} is proved in the same way as Proposition 3.4 of
\cite{LPP}. Let $\zeta=\zeta(p,\beta)=\PP^\beta(X_i\geq
1\, \forall i\geq 1)$.

\begin{cor}\label{regen}
The law of $\{(F_{R_n}),S_{Z_{R_n}}\}$ is absolutely continuous with respect to $\PP^\beta$. Furthermore, its
Radon--Nikodym derivative with respect to $\PP^\beta$ is
\begin{equation*}
\frac{d\mu}{d\PP^\beta}= {\bf I}_{\{X_i\geq
1\, \forall i\geq 1\}} \cdot\zeta^{-1}
\leq \zeta^{-1} < \infty\,.
\end{equation*}
\end{cor}

\section{Proof of Proposition \ref{mainlargep}}\label{proflrgp}
Proposition \ref{mainlargep} is a consequence of the following lemma.
\begin{lemma}\label{finexpec2} Let $\beta$ and $p$ be as in Lemma
\ref{sqrt}. Then, $E^\beta(R_{2}-R_1)<\infty$.
\end{lemma}
We will first show the following.
\begin{lemma}\label{finexpec1} Let $\beta$ and $p$ be as in Lemma
\ref{sqrt}. Then, $E^\beta (R_1)<\infty$.
\end{lemma}

\begin{proof}
We will show that
\begin{equation}\label{renbefn}
\sum_{n=1}^\infty \PP^\beta(R_1>n)<\infty.
\end{equation}

We will estimate $\PP^\beta(R_1>n)$ for $n$ large enough in order to
show (\ref{renbefn}).
Let $u$ be as in the proof of Lemma \ref{sqrt}.
Let $A_1$ be the event that the even
traps (as defined in Page \pageref{pageeventrap}) in $[-n,n]^2$ are of
size not larger than $ u\log n$. 
For $n$ large enough, the probability of $A_1$ is at least $1-n^{-2}$.

Let
\begin{equation*}
\kappa=K^2\max\left(\frac{3}{\log\beta},2u\right),
\end{equation*}
where $K$ is the constant from Lemma \ref{conduc}. Let $\gamma_n$ be 
the smallest even integer $\geq \kappa(\log n)^2$ and let
\begin{equation*}
T_i=\inf\{k:X_k\geq i\gamma_n \}
\end{equation*}
\begin{claim}\label{almostindep}
Let $g(z)$ be ${\bf 1}_{G(z)}$. Let $\eta =\eta(p)>0$ be the
probability of a vertex to be good. Then, there exists
$\eta'>\eta/2$ and i.i.d. Bernoulli random variables $L_i, i=1, 2, \ldots$
where $L_i=1$
with probability $\eta'$ and $L_i=0$ 
with probability $1-\eta'$,
such that the total variation
distance between the conditional distribution of $(g(Z_{T_i}),1\leq i \leq n^{1/20})$, given $A_1$,
and the distribution of $(L_i,1\leq i \leq n^{1/20}) $ is bounded
by $n^{-2}$.
\end{claim}
\begin{proof}
If we condition on nonexistence of even traps of 
size larger than $u\log n$, the
point $Z_{T_i}$ is good if and only if there exists a good path
starting at $Z_{T_i}$ and ending at the line
\[
\{((i+1)\kappa(\log n)^2, y) : y\in\sublattice\}.
\]
Let $\eta'$ be the $\PP^\beta$--probability of the
existence of such a path. 
We now define the random variables $L_i, i=1, 2, \ldots$:
let $L_i$ be the indicator of the event that there is a good path starting at
$Z_{T_i}$ and ending at the line $\{((i+1)\kappa(\log n)^2,y) :
y\in\sublattice\}$. Since we
condition on the origin being in the infinite cluster in the left half--plane,
the conditional probability of 
$\{L_n=1\}$, given ${\cal H}(Z_{T_n})$ and
the percolation configuration on $\{(x,y)| x\leq Z_{T_n} \}$, does not depend 
on ${\cal H}(Z_{T_n})$ and
the percolation configuration on $\{(x,y)| x\leq Z_{T_n} \}$. Therefore, the 
random variables $L_i$ 
are i.i.d.
Since the conditional distribution of 
$(g(Z_{T_i}),1\leq i \leq n^{1/20})$, given $A_1$,
was obtained from the distribution of $(L_i,1\leq i \leq n^{1/20})$
by conditioning on an event of
probability at least $1-n^{-2}$, the total variation
distance  between the two distributions  
is bounded
by $n^{-2}$.
\end{proof}

Let $A_2$ be the event that $T_i<n$ for every $1\leq i\leq n^{1/20}$.
By Lemma \ref{sqrt}, for $n$ large enough, $\PP^\beta(A_2)\geq 1-n^{-2}$.

Let $A_3$ be the event that there are at least $\frac{1}{4}\eta n^{1/20}$
values of $i$ in $\{1,2,...,[n^{1/20}\, ]\}$ such that $g(Z_{T_i})=1$.
By Claim \ref{almostindep}, for $n$ large enough, 
$\PP^\beta(A_3|A_1)\geq 1-3n^{-2}$ and therefore 
$\PP^\beta(A_3) \geq 1-4n^{-2}$.

Let $\xi_j$  be the $j$--th value of
$T_i$ such that $g(Z_{T_i})=1$. We define
\begin{equation*}
D(i)=\inf\{X_k-X_{\xi_i}:k>\xi_i\}
\end{equation*}
and 
\begin{equation*}
\tildd(i)=\inf\{X_k-X_{\xi_i}:\xi_i<k<\xi_{i+1}\}\,.
\end{equation*}

\begin{claim}\label{eprime}
There exists $\rho>0$ such that
\begin{equation*}
\PP^\beta_{{\cal H}(\xi_i)}(D(i)=1)\geq\rho \quad \PP^\beta \mbox{--a.s.}
\end{equation*}
\end{claim}
\begin{proof}
The claim is a direct consequence of Lemma \ref{conduc}: take $\ell$
such that $K\beta^{-\sqrt{\ell}/K} < 1$, then the probability of the event 
$\{D(i)=1\}$ is bounded below by 
$(\beta+3)^{-2\ell}(1-K\beta^{-\sqrt{\ell}/K})$.
\end{proof}

Obviously, $D(i)\leq\tildd(i)$ for every $i$. Therefore,
\begin{equation}\label{manju}
\PP^\beta_{\H(\xi_i)}(\tildd(i)=1)\geq\rho \quad \PP^\beta \mbox{--a.s.}
\end{equation}
Note that, for all $i$,
$\tildd(i-1)$ is $\H(\xi_{i})$--measurable. Therefore, by (\ref{manju}) and successive conditioning,
for every $k$,
\begin{equation*}
\PP^\beta(\tildd(i)<1 \text{ for all } i=1,\ldots,k) \leq (1-\rho)^k.
\end{equation*}

Let $A_4$ be the event that there exists some 
$1\leq i\leq \frac{1}{4}\eta n^{1/20} $
such that $\tildd(i)=1$.
Then,
\begin{equation*}
\PP^\beta (A_4) \geq 1-(1-\rho)^{\frac{\eta n^{1/20}}{4}}=1-o(n^{-2})\, .
\end{equation*}
Let $A_5$ be
the event that $D(i)=\tildd(i)$ for every $1\leq i\leq
\frac{1}{4}\eta n^{1/20}$.
For every $i$, 
\begin{equation}\label{min}
D(i)=\min\left(\tildd(i),D(i+1)+X_{\xi_{i+1}}-X_{\xi_i}\right)\, .
\end{equation}
By Lemma \ref{conduc}, 
\begin{equation}\label{fromcond}
\PP^\beta (D(i+1)\leq X_{\xi_i}-X_{\xi_{i+1}})\leq\PP^\beta(D(i+1)\leq-\kappa\log n)
\leq K\beta^{-3\frac{\log n}{\log\beta}} = Kn^{-3}
\end{equation}
Combining (\ref{min}) and (\ref{fromcond}), we get that
$\PP^\beta(D(i)\neq\tildd(i))\leq Kn^{-3}$ for every $i$. Therefore, 
$\PP^\beta(A_5)= 1-o(n^{-2})$.

\begin{claim}
If $A_1$, $A_2$, $A_3$, $A_4$ and $A_5$ all occur, then $R_1\leq n$.
\end{claim}
\begin{proof}
By the occurrence of $A_2$ and $A_3$, $\xi_i<n$ for every
$i\in B_n=[1,\ldots,\frac{1}{4}\eta n^{1/20}]$. By the occurrence of $A_4$, there
exists $i_0\in B_n$ such that $\tildd(i_0)=1$. By the occurrence of $A_5$,
$D(i_0)=1$. Let $t=\xi_{i_0}$. Then $t<n$. By the definition of
$\{\xi_i\}$, the epoch $t$ is a fresh epoch. On the other hand, for
every $k>t$,
\[
X_k\geq X_t+\min(X_j-X_t:  j>t)=X_t+D(i_0)=X_t+1>X_t
\]
and therefore $t$ is a regeneration epoch.
\end{proof}

Hence 
\begin{equation*}
\PP^\beta(R_1>n)\leq \PP^\beta(A_1^c)+\PP^\beta(A_2^c)+\PP^\beta(A_3^c)+
\PP^\beta(A_4^c)+\PP^\beta(A_5^c)=O(n^{-2}),
\end{equation*}
which yields (\ref{renbefn}).
\end{proof}
\begin{proof}[Proof of Lemma \ref{finexpec2}]
For a random variable $X$ and a distribution $\nu$, we denote the
expected value of $X$ under $\nu$ by $E_\nu(X)$.
We want to show that $E_{\PP^\beta}(R_2-R_1)<\infty$. Recall the
distribution $\mu= \mu^\beta$ from Section \ref{regenerations}. 
The distribution
of $R_2-R_1$ under $\PP^\beta$ is the same as the distribution  of $R_1$
under  $\mu$. Therefore, all we need to show is that
$E_\mu(R_1)<\infty$. But, using Corollary \ref{regen} and 
Lemma \ref{finexpec1},
\begin{equation*}
E_\mu(R_1)\leq E_{\PP^\beta}(R_1)\cdot
\sup\left(\frac{d\mu}{d\PP^\beta}\right)<\infty\, .
\end{equation*}
\end{proof}

\section{Renormalization}\label{renorm}
In this section we show how to combine standard renormalization ideas
with our arguments in order to carry over our results for every $p>p_c$.
We use the renormalization scheme that is used in \cite{BPP},
\cite{EO} and \cite{ABBP}. 
Fix a value $p\in(p_c,1)$.

Notice that everything we did so far is also valid when we consider site
percolation with retention probability $\hat{p}<1$ instead of bond percolation. 
We will assume that $\hat{p}<1$ is close enough 
to $1$ to apply our previous arguments (to be specified later).
Let $N$ be a
(large) positive integer, divisible by $8$. 

For each $v\in N\cdot\lattice$ define $Q_N(v)$ to be the square of side--length 
$5N/4$
centered at $v$. Let $p\in(p_c,1)$. Consider i.i.d. bond percolation
with parameter $p$ on $\lattice$, and let $A_p$ be the random
set of vertices $v\in N\cdot \lattice$ such that $Q_N(v)$ contains a connected
open component which connects all $4$ faces of $Q_N(v)$ but contains no other
connected open component of diameter greater than $N/10$. It follows from
Proposition 2.1 in 
\cite{AnPi} that if $N$ is large enough
then $A_p$ dominates i.i.d. site percolation with parameter $\hat{p}$ on
$\lattice$. We choose $N$ to be such a large enough value. 

For $p$ close to $p_c$ it is possible to show that there is (a.s.) no point in the lattice that
satisfies the definition of a good point (Definition \ref{def:goodpnt}
on page \pageref{def:goodpnt}).
Therefore, we need a new notion of a point being
good. In order to avoid confusion, we will use the term $p$--good for
the new definition.
\begin{defn}
We say that a square $Q_N(v)$ is {\bf $p$--good} if $v\in A_p$ and there
exist
$v_1=v=(x_1,y_1),v_2=(x_2,y_2),v_3=(x_3,y_3),v_4=(x_4,y_4),...$  such that
\\(A) For every $k$, $x_k-x_{k-1}=N$ and $y_k-y_{k-1}=\pm N$. 
\\(B) For every $k$, both $v_k$ and $v_k  + (N,0)$ are in $A_p$.

A square is considered {\bf $p$--bad} if it is not $p$--good.
\end{defn}
If $z$ is a point in the $p$--good square $Q_N(v_1)$ that belongs to the big component
in the square, then there exists an infinite path starting at $z$ that
is contained in the union of the squares $Q_N(v_1), Q_N(v_1+(N,0)),
Q_N(v_2), Q_N(v_2+(N,0)), Q_N(v_3), Q_N(v_3+(N,0)), \ldots$. (This follows from the definition of $A_p$ --- note that a connected component crossing the 
overlapping part of two good squares has to cross both squares!) We call
this path a {\bf $p$--good path} starting at $z$.
\begin{defn}
We say that a point $z$ is {\bf $p$--bad} if it is in a $p$--bad square and belongs to
the infinite cluster.
\end{defn}

\begin{defn}
We say that a point $z$ is {\bf $p$--good} if
\\(A) $z$ is not $p$--bad.
\\(B) There exists an (infinite) $p$--good path
$z_1=(x_1,y_1)=z,z_2=(x_2,y_2),z_3 =(x_3,y_3),z_4=(x_4,y_4),\ldots$
starting at $z$ such that $x_k>x_1$ for every $k>1$.
\end{defn}

\begin{defn}
A $p$--trap is a connected component of $p$--bad points.
\end{defn}

\begin{remark}
The reader is advised to notice
that:
\\(A) The squares are not disjoint. Therefore a point could belong to
both a $p$--good square and a $p$--bad square. In this case, if it is connected
to infinity then it is $p$--bad.
\\(B) Not all of the points that are connected to infinity are
$p$--good or $p$--bad.
\\(C) A $p$--good path may also contain $p$--bad points.
\\(D) If $\hat{p}$ is close enough to $1$, a square has a positive probability of being $p$--good, and a vertex
has a positive probability of being $p$--good. 
\end{remark}

In particular, a point at the boundary of a $p$--trap might not be a $p$--good 
point. Therefore, we also need the following weaker definition.
\begin{defn}
A {\bf $p$--OK} point is a point that is not $p$--bad and is in the big cluster of a $p$--good square.
\end{defn}

Once we defined a $p$--trap and a $p$--OK point, the argument for transience of the random walk follows
the same lines as in the case where $p$ is close enough to $1$. 
More precisely, let $T_p(z)$ be the $p$--trap containing $z$, and let $L_p(z)$  and
$W_p(z)$ be the length and the width of $T_p(z)$.

\begin{lemma}\label{smalltraps}
There exists $\alpha<1$ such that $\lP(L_p(0)\geq n)\leq\alpha^n$ and 
$\lP(W_p(0)\geq n)\leq \alpha^n$ for
every $n$.
\end{lemma}
The proof is the same as the proof of Lemma \ref{trapsize}, 
assuming $\hat{p}$ is close enough to $1$  and considering 
oriented percolation on the sublattice of the centers of squares.

\begin{lemma}\label{goback}
For a point $z=(x,y)$, let $\OK(z)$ be the event that $z$ is a $p$--OK
point. Then,
there exists a constant $K'=K'(p,\beta)$
such that for every $\ell=1,2,\ldots$,
\begin{equation*}
\PP^\beta_{\H(n)}\left(\hbox {\rm there is an } m\geq n 
\hbox{ \rm such that } X_m\leq
x-\ell)\left|Z_n=z \right., F(n),\OK(z)\right)\leq
K'\beta^{-\sqrt{\ell}/K'}\  .
\end{equation*}
\end{lemma}
The proof, again, is similar to that of Lemma \ref{conduc} since one
can bound from below the conductance of every $p$--good path starting
at $z$.

In order to prove the equivalents of Lemma \ref{sqrt} and Lemma
\ref{transtoright} we need the following simple claim:

\begin{claim}\label{OK}
Let $T$ be a $p$--trap. Every point at the boundary of $T$ is
$p$--OK.
\end{claim}

Using Claim \ref{OK}, Lemma \ref{goback} and Lemma \ref{smalltraps}
we can now prove
the following two lemmas the
same way Lemma \ref{sqrt} and Lemma \ref{transtoright} were proved.

\begin{lemma}\label{rensqrt}
For $\beta >1$  close enough to
$1$, there exists a constant $C$ such that for every $n$ large enough,
\begin{equation*}
\PP^\beta(X_n\leq Cn^\frac{1}{10})\leq n^{-2}.
\end{equation*}
\end{lemma}

\begin{lemma}\label{rentranstoright}
Let $\beta>1$, then
\begin{equation}
\lim_{n\to\infty}X_n=\infty \quad \PP^\beta\mbox{--a.s.}
\end{equation}
\end{lemma}

The proof of Theorem \ref{main} now a follows the same lines as the proof of
Proposition \ref{mainlargep}
in Section \ref{proflrgp}, using the notions ``$p$--good'' and
``$p$--trap'' instead of ``good'' and ``trap''.

\section{Zero speed region}\label{speedzero}
\begin{thm}\label{zerospeed}
For every $ p \in (p_c, 1)$,  
there exists a finite value $\beta_u = \beta _u(p) > 1$ such that for 
$\beta > \beta_u$,
\begin{equation*}
\lim_{n\to\infty}\frac{X_n}{n}= 0\quad \PP^\beta\mbox{--a.s.}
\end{equation*}
Further, 
$\lim_{p \searrow p_c} \beta_u(p) =1$.
\end{thm}
\begin{proof}
We will first define $\beta_u$ and show that for $\beta > \beta_u$, 
the speed of the random walk is $0$.
For this purpose, we will consider configurations where the origin $0$ is 
the beginning of a dead end. Call a vertex $z=(x,y)$ the {\bf beginning of a dead end} if $z$ is in the 
infinite cluster to its left, but in a finite cluster to its
right. The dead end starting at $z$ is the finite cluster to the
right of $z$, containing $z$. We now consider a dead end $A$ starting at the origin. 
Let $d(A): = \max\{x: (x,y)\in A\}$ denote the depth 
of $A$. Let $N(A)$ denote the number of vertices of $A$ which are on
the line $\linezero=\{(0,y):y\in\sublattice\}$.
Let $E_A$ denote the set of edges of $A$ and $B_A$
denote the set of all edges which have at least one vertex in $A$, 
but are not in $E_A$.
The probability of $A$ under i.i.d. bond percolation is
$p_A := p^{|E_A|}(1-p)^{|B_A|}$. Let $DE$ denote the set of all dead ends.
The following claim is easy, we omit its proof. 
\begin{claim}\label{dcinf}
Let $\{\omega_r = A\}$ denote the event that all the edges in $E_A$ are 
open in $\omega$ and all the edges in $B_A$ are closed in $\omega$. Then,
$$
\lP\left(\omega_r=A\right)\geq C(p) p_A
$$
where $C(p)$ is a constant depending only on $p$.
\end{claim}
Let 
\begin{equation}\label{Gammadef}
\Gamma(p, \beta) : = \sum\limits_{A \in DE}p_A 
\left(N(A)^{-1}
\sum\limits_{e=(z_1, z_2) \in E_A}
\beta^{x_1 \vee x_2}
\right)
\end{equation}
where $z_1= (x_1,y_1)$ and $z_2= (x_2,y_2)$. We define $\beta_u = \beta_u(p)$
 as the threshold value for convergence, i.e. such that 
$\Gamma(p, \beta) < \infty$ for 
$\beta < \beta_u$ and $\Gamma(p, \beta) = \infty$ for $\beta > \beta_u$.
It is easy to see, giving a lower bound for $\Gamma(p, \beta)$, that 
$\beta_u< \infty$ for all $p$. Let $T_0: =\inf\{j>1: X_j =0\}$, and let $T_A$ be the time spent in the 
dead end $A$. Then, on $\{\omega_r=A\}$, 
$E^\beta_\omega(T_A)= E^\beta_\omega(T_0| X_1\geq 0)$. 
\begin{lemma}\label{inftytime}
For $\beta > \beta_u$, 
\begin{equation}\label{inftime}
E^\beta (T_0) = \infty\,.
\end{equation}
\end{lemma}
\begin{proof}
We will show that for $\beta > \beta_u$, the expected time spent in a 
dead end starting at $0$ is infinite, 
giving a lower bound for the latter by considering the time spent in 
the dead end up to the first 
return to $\linezero$.
Consider the random walk on $A$, starting 
from $0$.
Let $T_{A,0}: =\inf\{j>1: Z_j \in\linezero\}$. 
We have, on $\{\omega_r=A\}$, 
\begin{equation}\label{timeinA}
E^\beta_\omega(T_{A,0}| X_1\geq 0) \geq
\frac{2}{3+ \beta}N(A)^{-1}\sum\limits_{e=(z_1, z_2) \in E_A} \beta^{x_1 \vee x_2} \,
\end{equation}
This follows from the fact that for a recurrent Markov chain on $A$ with 
invariant measure $\pi$, the expected return time to a vertex $z$ is
$\pi(A)/\pi(z)$.
In our case, the invariant measure
$\pi(z)$ is
given by the sum of the weights of all edges $e=(z, \cdot)$ where the
weight of an edge
$e=(z_1, z_2)$, $z_1= (x_1,y_1)$, $z_2= (x_2,y_2)$ is given by $\beta^{x_1 
\vee x_2}$, hence $\pi(A) = 2\sum\limits_{e=(z_1, z_2) \in E_A} \beta^{x_1 \vee x_2}$.
(\ref{timeinA}) now follows by merging all of the vertices of 
$A\cap \linezero $ into one vertex.
\end{proof}

\begin{lemma}\label{inftytimediff}
For $\beta > \beta_u$,  the speed of the random walk is zero.
\end{lemma}
\begin{proof}
We define a sequence of {\bf ladder times} $L_1, L_2, \ldots $.
Let $L_1$ be the first fresh epoch such that $Z_{L_1}$ is the beginning of 
a dead end. Let $A_{L_1}$ be the dead end starting at $L_1$ and $d(A_{L_1})$ its depth.
Let $L_2$ be the first fresh epoch such that $X_{L_2} > X_{L_1} + d(A_{L_1})$ 
and 
$Z_{L_2}$ is the beginning of 
a dead end, and continue the recursion. If $n$ is a fresh epoch,
the environment to the right of $Z_n$ has the same distribution as the environment to the right of the 
origin under $\lP$. 
Therefore, the probability that
the first hitting time of 
$\{(x,y): x= X_{L_i} + d(A_{L_i})+1\}$ is a ladder 
time is
strictly positive and does not depend on $i$.
In particular, there are infinitely many ladder times.
We will show that $X_{L_n}/L_n \to 0$, $\PP^\beta$--a.s. for $n \to \infty$.
Note that $L_{i+1} -L_i \geq T_{A_i}$ and the random variables $(T_{A_i})$ are
i.i.d. under $\PP^\beta$  and have, due to Lemma \ref{inftytime}, infinite expectation for $\beta > \beta_u$.
This implies that $L_n/ n \to \infty$, $\PP^\beta$--a.s. for $n \to \infty$. 
On the other hand,
the random variables $X_{L_{i+1}} - X_{L_i}$ are i.i.d. and we
claim that they have exponential tails and, in particular, 
finite expectations. 
To see this, note that due to Lemma \ref{trapsize} and Lemma \ref{smalltraps},
the depth of a dead end has an exponential tail, 
i.e. $\lP(d(A_0) \geq s)\leq \exp(-c(p)s)$ for $s$
large enough, where $c(p)$ is some constant depending only on $p$.
For an integer $t$ which is divisible by $20$, we want to estimate the
probability of the event 
\[
X_{L_{i+1}}-X_{L_i}>t.
\]
Let $s = t/20$.
Let $\tau_j:=
\inf\{k: X_k = X_{L_{i}} + 10 j\}$, $j=1, 2, \ldots $.
Let $B$ denote the event that $0$ is connected to 
$\{(10, y): y \in \sublattice \}$
if we
remove all the vertices on the line $\{(-1, y): y \in \sublattice \}$, 
and let $\gamma= \lP(B)$. Then, conditioning on the event that the dead end beginning at $L_{i}$ has depth at most $\hetzi t$,
consider the fresh epochs $\tau_j, j= 11s, \ldots, 20s$. They have either to be beginnings of dead 
ends or they have to be connected to the next line at distance $10$. 
Hence
$$
\PP^\beta(X_{L_{i+1}} - X_{L_i} \geq  t)
\leq 
\exp\left(-c(p)\hetzi t\right) + \gamma^{s} \leq \exp(-\widetilde c(p)t)
$$
for some constant $\widetilde c(p)$.

Hence, $\limsup X_{L_n}/n < \infty$, $\PP^\beta$--a.s. and we conclude that
$X_{L_n}/L_n \to 0$, $\PP^\beta$--a.s.
Since $L_{n+1} /L_n \to 1$, $\PP^\beta$--a.s. for $n \to \infty$, this suffices to prove that $X_n/n \to 0$, $\PP^\beta$--a.s. for $n \to \infty$. 
\end{proof}

\begin{lemma}\label{ptocrit}
We have $\beta_u(p) \to 1$ for $p \searrow p_c$.
\end{lemma}
\begin{proof} Fix $\beta>1$.
Let 
$\partial_+ B_n : =\{(x,y): x=n\ve |y| \leq n\}$.
Then, for every $n$, using the proof of Lemma 
\ref{inftytime},
\begin{eqnarray*}
\label{split}
& & E^\beta(T_0| X_1 \geq 0)\\
&\geq&\beta^n \lP( 0 \hbox{ \rm is connected to a vertex  } v 
\in \partial_+ B_n)\\
& &\times \lP(\hbox{\rm the clusters of all vertices } z \in \partial_+ B_n 
\hbox{ \rm are finite})\, .
\end{eqnarray*}
Now, since $p > p_c$, 
\begin{equation}
\label{connex}
\lP( 0 \hbox{ \rm is connected to a vertex  } v \in \partial_+ B_n) \geq
\mu_p >0.
\end{equation}
Let $\delta >0$ be such that
\begin{equation*}
W:=\beta(1-\delta)^4>1 .
\end{equation*}
For $p$ close enough to $p_c$, since $\theta(p_c)=0$, 
$P_p(C_0 \hbox { \rm finite}) \geq 1-\delta$ (where $\theta(p)$
denotes the probability that the origin belongs to an infinite open
cluster,
and we refer to \cite{grimmett} for the fact that $\theta(p_c)=0$). Hence,
using the FKG inequality, $P_p(\hbox{\rm the clusters of all vertices } 
z \in \partial_+ B_n \hbox{ \rm are finite})$ 
can be estimated as follows. For $p$ close enough to $p_c$,
$$
P_p(\hbox{\rm the clusters of all vertices } z \in \partial_+ B_n 
\hbox{ \rm  are finite})
\geq (1-\delta)^{4n}\, .
$$
We conclude that also
$$
\lP(\hbox{\rm the clusters of all vertices } z \in \partial_+ B_n 
\hbox{ \rm  are finite})
\geq c(1-\delta)^{4n}\, .
$$
for some constant $c=c(p)$. Thus, for every $n$,
\begin{equation}\label{forn}
E^\beta(T_0| X_1 \geq 0) \geq \mu_p\beta^n(1-\delta)^{4n}=\mu_pW^n.
\end{equation}
Since $W>1$ and (\ref{forn}) holds for every $n$, we conclude that 
$E^\beta(T_0| X_1 \geq 0) = \infty$. Recalling (\ref{Gammadef}) and
(\ref{timeinA}), we see that $\Gamma(p, \beta) = \infty$, hence 
$\beta \geq \beta_u$.   
\end{proof}

Theorem \ref{zerospeed} now follows from Lemma \ref{inftytimediff} and
Lemma \ref{ptocrit}.
\end{proof}


\end{document}